\documentclass[12pt]{amsart} 
\usepackage{amsthm, amscd}
\usepackage{amssymb,latexsym} 

\theoremstyle{plain}
\newtheorem*{introtheorem}{Theorem}

\theoremstyle{definition}

\theoremstyle{remark}

\newcommand{\OO}{\mathcal{O}}

\begin{document}

\title{\bf On the base locus of the generalized theta divisor}
\author{Mihnea Popa}
\address{Department of Mathematics
\\University of Michigan \\ Ann Arbor, MI 48109 \hspace{.3cm} and} 
\address{Institute of
Mathematics of the Romanian Academy \\ Calea Grivitei 21 \\ Bucharest 
\\ Romania} 
\email{mpopa@math.lsa.umich.edu}

\maketitle

{\bf{Abstract.}}  In response to a question of Beauville, we give a     
new class of examples of base points for the linear system $|\Theta|$ on
the moduli space $SU_{X}(r)$ of semistable rank $r$ vector bundles of
trivial determinant on a curve $X$ and we prove that for sufficiently
large $r$ the base locus is positive dimensional.

\section{Introduction}

Let $X$ be a compact Riemann surface of genus $g$. In his survey
\cite{Beauville}, A. Beauville raises a few questions about the base locus
of  
the linear system $|\Theta|$, where $\Theta$ is the theta (or
determinant) bundle on the moduli space $SU_{X}(r)$ of semistable
rank $r$ vector bundles on $X$ of trivial determinant. It is known
(see \cite{Beauville}, \S3) that $E\in SU_{X}(r)$ is a base point for
$|\Theta|$ if
and only if $H^{0}(E\otimes L)\not=0$ for every $L\in$
Pic$^{g-1}(X) $\footnote{When $E$ is semistable, this should be
interpreted as a statement about its equivalence class: by a common
argument using the Jordan-H\"older filtration of $E$ and the fact
that $\chi(E\otimes L)=0$, it is easy to see that $H^{0}(E\otimes
L)\neq0, \forall L\in {\rm Pic}^{g-1}(X)$ iff $H^{0}(gr(E)\otimes 
L)\neq0, \forall L\in {\rm Pic}^{g-1}(X)$. So this property does not
depend on the choice of a bundle in the class of $E$.}.
When $r=2$, or $r=3$ and $X$ is of genus $2$ or generic of any genus,
it is known that $|\Theta|$ is base point free. However, M. Raynaud
constructs in
\cite{Raynaud} examples of bundles which lead to the existence of base
points of $|\Theta|$ for
$r= n^{g}$, where $n$ is an integer $\geq 2$ dividing $g$. His
construction gives finitely many base points in each genus.
Among other things, Beauville asks in \cite{Beauville} if one could find
new examples
of such base points and if the base locus is actually of strictly positive
dimension. 

The purpose of this note is to give at least a partial answer  
to Beauville's question. From a qualitative point of view, our results
can be summarized in the following:
\begin{introtheorem}
(a) For every $g\geq 2$, there exists a rank $\rho(g)$ such that for  
all $r\geq \rho(g)$ the linear system $|\Theta|$ on $SU_{X}(r)$ has 
base points. Also, for every $g\geq 2$ there exist ranks where
some base points are stable.
 
(b)  Moreover, for every $g\geq 2$ and every $k\geq 2$, there exists an
integer 
$\rho(k,g)$ such that for all $r\geq \rho(k,g)$, the base locus
of $|\Theta|$ on
$SU_{X}(r)$ has dimension at least $(k-1)g$.
\end{introtheorem}

For the examples and for more precise numerical statements see Section 2.

Another question raised by Beauville adresses the freeness of
$|2\Theta|$. In the spirit of \cite{Donagi}, one could also look at other
moduli
spaces, not
necessarily of trivial determinant. We remark, by using a theorem of Lange
and Mukai-Sakai (see \cite{Lange},\cite{Mukai}), that the global
generation of low
multiples of the theta divisor on such moduli spaces cannot go hand in
hand with the strange duality conjecture (see \cite{Beauville}, \S8).

\textbf{Acknowledgement:} I would like to thank Rob Lazarsfeld
for suggesting the problem. I would also like to thank him and Iustin 
Coand\u a for valuable discussions.

\section{The examples and a lower bound for the dimension of the base
locus}

Throughout the paper we will denote by
$SU_{X}(r,A)$ the moduli space of
semistable bundles of rank $r$ and fixed determinant $A$ and by
$U_{X}(r,d)$
the moduli space of semistable bundles of rank $r$ and degree $d$ on $X$.
 
Consider a line bundle $L$ on $X$ of degree $d\geq 2g+1$ . We
will restrict to the case $g\geq 2$, since for $g\leq 1$ the
space $SU_{X}(r)$ is well understood. Denote by $M_{L}$ the   
kernel of the evaluation map:

$$0 \longrightarrow M_{L} \longrightarrow H^{0}(L)\otimes \mathcal{O}_{X}
\longrightarrow L \longrightarrow 0$$
and let $Q_{L}= M_{L}^{\ast}$. These vector bundles are well known for
their importance in the study of the minimal resolution of $X$ in
the embedding defined by $L$ (see \cite{Lazarsfeld}, \S1 for a survey).

Among the properties of  $Q_{L}$, we quote from \cite{Lazarsfeld}, \S1.4,
the
following: if
$x_{1},\ldots,x_{d}$ are the points of a generic hyperplane section of
$X\subset {\mathbb P}(H^{0}(L))$, then $Q_{L}$ sits in an extension:
$$0 \longrightarrow 
\bigoplus_{i=1}^{d-g-1}\mathcal{O}_{X}(x_{i})
\longrightarrow
Q_{L} \longrightarrow \mathcal{O}_{X}(x_{d-g}+\ldots +x_{d})
\longrightarrow 0$$
This induces for every integer $p$ an inclusion:
\begin{equation} 
0 \longrightarrow \underset{1\leq i_{1}<\ldots
<i_{p}\leq
d-g-1}{\bigoplus}\mathcal{O}_{X}(x_{i_{1}}+\ldots +x_{i_{p}})
\longrightarrow
\bigwedge^{p}Q_{L}
\end{equation}
\newline 
Recall also from \cite{Ein}, \S3 that $Q_{L}$ is stable, and so
$\bigwedge^{p}Q_{L}$ is poly-stable (i.e. a direct sum of
stable bundles of the same slope).
\medskip
\newline
\textbf{Definition.} Similarly to a definition in \cite{Raynaud}, we
say that a vector bundle $E$
satisfies property $(\ast)$ if and only if:
$$H^{0}(E\otimes \xi) \neq  0 , \forall 
\xi\in {\rm  Pic}^{0}(X).$$
In all that follows we will denote
$\gamma := [\frac{g+1}{2}]$.
\begin{proof}[Proof of Theorem.(a).] Notice first that to find a base
point for
$|\Theta|$ 
it is enough to exhibit a semistable bundle $E$ of integral slope $0\leq
\mu(E)\leq g-1$
satisfying $(\ast)$, since we could then twist by a suitable line bundle.
\medskip
\newline
\textbf{Claim:} \emph{For every line bundle $L$ on $X$ of degree $d\geq
2g+1$, the bundle
$\bigwedge^{\gamma}Q_{L}$ satisfies property $(\ast)$}.
\medskip
\newline
\emph{Proof of claim.}
From (1) it is clear that for $x_{1},\ldots ,x_{p}$ general points on $X$
we have $H^{0}(\bigwedge^{p}Q_{L}(-x_{1}-\ldots$
$-x_{p}))\neq0$. 
So for any $p$, a generic line bundle $\xi\in$ Pic$^{0}(X)$ of the form
$\xi =
\OO_{X}(A_{p}- B_{p})$, with $A_{p},B_{p}$ generic effective divisors of
degree $p$, satisfies
$H^{0}(\bigwedge_{}^{p}Q_{L}\otimes \xi)\neq 0$.
 
On the other hand it is well known (see \cite{Harris}) that every
$\xi\in$
Pic$^{0}(X)$ can be written in the form $\xi = \OO_{X}(A_{\gamma}-
B_{\gamma})$ with $A_{\gamma},B_{\gamma}$ effective divisors of degree
$\gamma$.
Hence
$H^{0}\left(\bigwedge^{\gamma}Q_{L}\otimes
\xi\right)\neq 0$ for a general $\xi$ and by semicontinuity the
same must hold for every  $\xi\in$ Pic$^{0}(X)$, which proves the claim.
\medskip

So, as noted above,
it is enough to get integral slopes for $\bigwedge^{\gamma}Q_{L}$ for
suitable choices of $d$. Since the
computations differ from case to case and tend to get messy, we will
restrict to giving examples that work uniformly rather then trying to find
the smallest possible rank for each genus.

We can obtain an uniform answer by choosing $d=
g\left(\gamma +1\right)$,
when we get
$\mu \left(\bigwedge^{\gamma}Q_{L}\right)= \gamma +1$.
The
corresponding rank will
be $rk\left(\bigwedge^{\gamma}Q_{L}\right)=
\left({g\cdot\gamma} \atop {\gamma}\right)$ (actually for most $g$'s 
this is by no means the best answer).

It is easy to see that since the bundles $\bigwedge^{\gamma}Q_{L}$ are   
poly-stable and satisfy $(\ast)$, at least one of their stable summands
(which have the same slope) must also satisfy $(\ast)$. Thus the
constructions above also give us examples of
stable base points in each genus. On the other hand, the existence of a
base point $E\in SU_{X}(r)$ induces
the
existence of decomposable base points for every rank $r^{'}\geq r$
:
simply take $E\oplus \mathcal{O}_{X}^{\oplus (r^{'}-r)}$.
\end{proof}
\medskip
{\bf{Remark}}: There are many versions of this construction that give
additional examples.
Let us just mention them without getting into numerology. One could
look at $\bigwedge^{p}Q_{L}$ for $p>\gamma$ such that
$\mu\left(\bigwedge^{p}Q_{L}\right)\leq g-1$ or work
with $S^{p}Q_{L}$ instead of $\bigwedge^{p}Q_{L}$. It is probably most
interesting though to replace $Q_{L}$ by $Q_{E}$, where $E$ is a
semistable bundle
of slope $\mu(E)>2g$ (so automatically very ample) and $Q_{E}$ is defined
exactly as $Q_{L}$. By \cite{Butler} $Q_{E}$ is known to be semistable and
a closer
analysis
shows that a result analogous to the claim above holds. Using this
construction one can check that by good numerical choices we can make
$\bigwedge^{\gamma}Q_{E}$ have any integral slope
$[\frac{g+1}{2}]< \mu \leq g-1$.
\medskip

The additional feature that makes these examples interesting is that they
come in positive dimensional families (roughly speaking by varying $L$),
so in the range covered by them the base locus is indeed positive
dimensional. 
\begin{proof}[Proof of Theorem.(b).] The following is the more precise
statement referred to in the introduction:
\newline
\textbf{Claim:} 
\emph{Fix $g\geq 2$, $d\geq 2g+1$, $k\geq 2$ and let $L$ be any line
bundles of degree $kd$. Then there exists a
$(k-1)g$ dimensional family of
(equivalence classes of ) semistable
bundles of rank $\left({k(d-g)} \atop {\gamma}\right)$ and fixed
determinant
$ L^{\otimes\left({k(d-g)-1} \atop
{\gamma -1}\right)}$ satisfying property $(\ast)$}.
\medskip
\newline
\emph{Proof of claim.}
Fix $L$ of degree $kd$. To every $k-1$ line bundles $L_{1},\ldots
,L_{k-1}\in {\rm Pic}^{d}(X)$ associate $L_{k}:= L\otimes
L_{1}^{\ast}\otimes\ldots\otimes
L_{k-1}^{\ast} \in {\rm Pic}^{d}(X)$, so that $L_{1}\otimes\ldots\otimes
L_{k}=
L$. Set $F_{L_{1},\ldots ,L_{k-1}}:= L_{1}\oplus\ldots\oplus L_{k}$.
Thus ${\rm det}(F_{L_{1},\ldots ,L_{k-1}})= L$ and clearly
$Q_{F_{L_{1},\ldots
,L_{k-1}}}= Q_{L_{1}}\oplus\ldots\oplus Q_{L_{k}}$ (cf. the remark above 
for the definition).

It is enough to prove that the morphism:
$$
\begin{array} {ccc}
\psi :~ {\rm Pic}^{d}(X)\times\ldots\times {\rm Pic}^{d}(X) &
\longrightarrow &
SU_{X}\left(\left({k(d-g)} \atop {\gamma}\right), L^{\otimes
\left({k(d-g)-1} \atop
{\gamma -1}\right)}\right) \\
 \!\!\! \\
{\left( L_{1},\ldots ,L_{k-1}\right)} & \rightsquigarrow &
\bigwedge^{\gamma}Q_{F_{L_{1},\ldots ,L_{k-1}}}
\end{array}
$$
is finite. We have:
$$
\bigwedge^{\gamma}Q_{F_{L_{1},\ldots ,L_{k-1}}}\cong 
\underset{i_{1}+\ldots +i_{k}=
\gamma}{\bigoplus}\left(\bigwedge^{i_{1}}Q_{L_{1}}\otimes\ldots\otimes 
\bigwedge^{i_{k}}Q_{L_{k}}\right)
$$
In particular $\bigwedge^{\gamma}Q_{F_{L_{1},\ldots ,L_{k-1}}}$
is poly-stable of slope $\gamma\cdot{\frac{d}{d-g}}$.

Assume now that:
$$
\bigwedge^{\gamma}Q_{F_{L_{1},\ldots ,L_{k-1}}}\cong
\bigwedge^{\gamma}Q_{F_{L_{1}^{\prime},\ldots ,L_{k-1}^{\prime}}} 
$$
for some other $L_{1}^{\prime},\ldots ,L_{k}^{\prime}$ as before.
By the previous formula one has inclusions:
$$
\bigwedge^{\gamma}Q_{L_{i}^{\prime}}\hookrightarrow
\bigwedge^{\gamma}Q_{F_{L_{1}^{\prime},\ldots
,L_{k-1}^{\prime}}}.
$$
As noted before, all the bundles above are
poly-stable (of the same slope ), so
$\bigwedge^{\gamma}Q_{L_{i}^{\prime}}$ is a direct sum of some collection
of the stable summands of
$\bigwedge^{\gamma}Q_{F_{L_{1},\ldots ,L_{k-1}}}$. There are finitely many
ways in which this can occur, so it is enough then to notice that the
morphism:
$$
\begin{array} {ccc}
\phi :~ {\rm Pic}^{d}(X)
 & \longrightarrow &
 U_{X}\left( \left({d-g} \atop
{\gamma}\right) ,
\gamma\cdot\frac{d}{d-g}\cdot\left({d-g} \atop
{\gamma}\right)\right) \\
M & \rightsquigarrow &  
\bigwedge^{\gamma}Q_{M}
\end{array}
$$
is finite. This is clear since det$(Q_{M})$=$M$ and the claim is proved.
\medskip

Taking in particular $d= g(\gamma +1)$ as in part (a), we get that the
bundles $\bigwedge^{\gamma}Q_{F_{L_{1},\ldots ,L_{k-1}}}$ have integral
slope $\gamma +1$ and so they lead to base points as before. Hence we can
take $\rho(k,g)= \left({k\cdot g\cdot \gamma} \atop {\gamma}\right)$.

This argument actually gives a statement about equivalence classes, since
the bundles in the family that we have constructed are all poly-stable.
Again, by adding trivial bundles we get the same statement in all ranks
$r\geq \rho(k,g)$. 
\end{proof}

One could conjecture that the base locus is at least $g$-dimensional
whenever it is non empty. Perhaps in view of the remark above an even 
more optimistic guess could be made. 

\section{Strange duality versus freenes of low multiples of theta}

In connection with Beauville's question about low multiples of $\Theta$,
we show that the strange duality conjecture implies the existence of base 
points on $|k\Theta|$ for small $k$, on suitable moduli spaces. 

Consider first, in general, the moduli space $SU_{X}(r,A)$
for some $A\in$ Pic$^{m}(X)$ , $m\in \mathbb{N}$, $m\leq g-1$. Consider
also
$F\in
U_{X}(k,k(g-1-m))$ and define $\Theta_{F}$ on $SU_{X}(r,A)$ to be
$\Theta_{F}=
\tau_{F}^{\ast}\Theta$, where $\tau_{F}$ is the map :
 
$$\begin{array}{ccc}\tau_{F} : SU_{X}(r,A)& \longrightarrow &
U_{X}(kr,kr(g-1))\\
 E~ &\rightsquigarrow & E\otimes F \end{array}$$
and $\Theta$ is the canonical theta divisor on $U_{X}(kr,kr(g-1))$. Set 
theoretically of course  $\Theta_{F}=\{ E\hspace{0.1cm}|\hspace{0.1cm} 
H^{0}(E\otimes
F)\not= 0 \}$. The
famous strange duality conjecture, discussed at length in \cite{Donagi},
or more
precisely its geometric formulation (see~\cite{Beauville}, \S8), asserts
that the
linear
system $|k\Theta|$ on $SU_{X}(r,A)$ is spanned by the divisors
$\Theta_{F}$ as $F$ varies in $U_{X}(k,k(g-1-m))$.

Let us consider in particular $E$ to be one of Raynaud's examples (see
\cite{Raynaud}, \S3)
i.e. a bundle $E$ with $rk(E)=n^{g}$, $\mu(E)=\frac{g}{n}$, $n|g$ which
satisfies $(\ast)$. A theorem of Lange and
Mukai-Sakai (see~\cite{Lange}, \cite{Mukai}) implies that every $F\in
U_{X}(k,k(g-1-\frac{g}{n}))$ has a subbundle $M\hookrightarrow F$ of
degree deg$(M)\geq \frac{k(g-1-\frac{g}{n})-(k-1)g}{k}$ . Assume further
that $\frac{g}{k}\geq 1+\frac{g}{n}$, which can be achieved for good
choices of $g$ and $n$. Then $m :=$deg$(M)\geq 0$ and so, by property
$(\ast)$: $$H^{0}(E\otimes M) = H^{0}(E\otimes M(-m{\rm p})\otimes
O_{X}(m{\rm p}))\not= 0 ,$$ where p is a point on $X$.
Since $M\hookrightarrow F$, we obtain:
$$H^{0}(E\otimes F)\not= 0 ~{\rm for~all}~ F\in
U_{X}\left(k,k\left(g-1-\frac{g}{n}\right)\right).$$
By the discussion above, the strange duality conjecture then implies
that $|k\Theta|$ on \\ $SU_{X}(n^{g},{\rm det}(E))$ has a base point at
$E$.
\medskip
\newline 
{\bf{Remarks}}: 1. The conclusion above suggests (assuming  
the strange duality conjecture is true!) that one should expect
$|k\Theta|$ to have base
points, say  for example for $k$ small enough with respect to $g$ or $r$,
even extrapolating to $SU_{X}(r)$.
\newline
2. In the discussion above we cannot use the examples from the
previous section instead
of Raynaud's examples, since the condition deg$(M)\geq 0$ is not
necessarily satisfied any more. 
\medskip
\newline
Let us conclude with another analogous application of the strange duality
conjecture:

Consider $F_{L}= \bigwedge^{\gamma}Q_{L}$  from the previous
section (of integral slope) and denote
$A_{L}={\rm det}(F_{L})= L^{\otimes \left({d-g-1} \atop
{\gamma -1}\right)}$, $r^{\prime}=
rk(F_{L})= \left({d-g} \atop {\gamma}\right)$. Then exactly by the
same argument as above, one can check that $F_{L}$ is a base point for
$|\Theta|$ on $SU_{X}(r^{\prime},A_{L})$ under mild assumptions on $d$.
One can also apply the same argument for at least the Raynaud examples
such that $g\geq n$ and $(g,n)\neq1$ (in particular for those of integral
slope).

\end{document}